\documentclass[a4paper,12pt]{amsart}
\usepackage{amsmath,amsfonts,latexsym,graphicx,amssymb,url}
\usepackage{amsmath,txfonts,pifont,pxfonts,manfnt}
\usepackage{psfrag}
\usepackage{wrapfig, subfigure,graphicx}
\usepackage{hyperref}
\usepackage[active]{srcltx}
\usepackage{pdfsync}
\usepackage{fixmath}
\usepackage{xcolor}
\newcommand{\R}{{\mathbb R}} 

\renewcommand{\(}{\left(}
\renewcommand{\)}{\right)}

\theoremstyle{plain}
\newtheorem{theorem}{Theorem}[section]

\newtheorem{remark}[theorem]{Remark}

\begin{document}

\setcounter{equation}{0}

\title[Reflection and refraction problems for metasurfaces]{Reflection and refraction problems for metasurfaces related to Monge-Amp\`ere equations}
\author{Cristian E. Guti\'errez and Luca Pallucchini}
\thanks{Research supported by NSF grant DMS--1600578}.
\address{Department of Mathematics, Temple University, Philadelphia, PA, 19122}
\email{gutierre@temple.edu, luca.pallucchini@temple.edu}

\begin{abstract}
A metasurface is a surface, tipically a plane, on which a function called phase discontinuity is chosen so that the metasurface produces a desired reflection or refraction job. 
We derive the equations that the phase discontinuity function must satisfy, they are Monge-Amp\`ere partial differential equations, and we prove existence of solutions.


\end{abstract}

\maketitle
\tableofcontents
\section{Introduction}

The subject of metalenses is a flourishing area of research having multiple potential applications and is one of the nine runners-up for Science's Breakthrough of the Year 2016 \cite{science-runner-ups-2016}. 
A typical problem is given a surface, in most applications a plane, find a function on the surface (a phase discontinuity) so that the pair, surface together with the phase discontinuity (the metalens) refracts or reflects light in a desired manner. 
This leads to ultra thin optical components that produce abrupt changes over the scale of the free-space wavelength in the phase. 
This is in contrast with classical lens design, where the question is to engineer the gradual accumulation of phase delay as the wave propagates in the device, reshaping the scattered wave front and beam profile at will. In particular, in lenses light propagates over distances much larger than the wavelength to shape wavefronts.
A vast literature appeared in recent years on this subject, see for example \cite{yu2011light}, \cite{aieta2012out},  \cite{aieta2012reflection}, and the comprehensive review articles \cite{yu-capasso:flatoptics} and \cite{chen2016review}.
For more recent work in the area and an extensive up to date bibliography, we refer to \cite{2107planaroptics:capasso}; see also  \cite{2016metalensesvisiblewavelenghts:capasso}-\cite{She_adaptive} and \cite{Chen_broadband}.

\begin{figure}[htbp]
\centering
\includegraphics[width=3in]{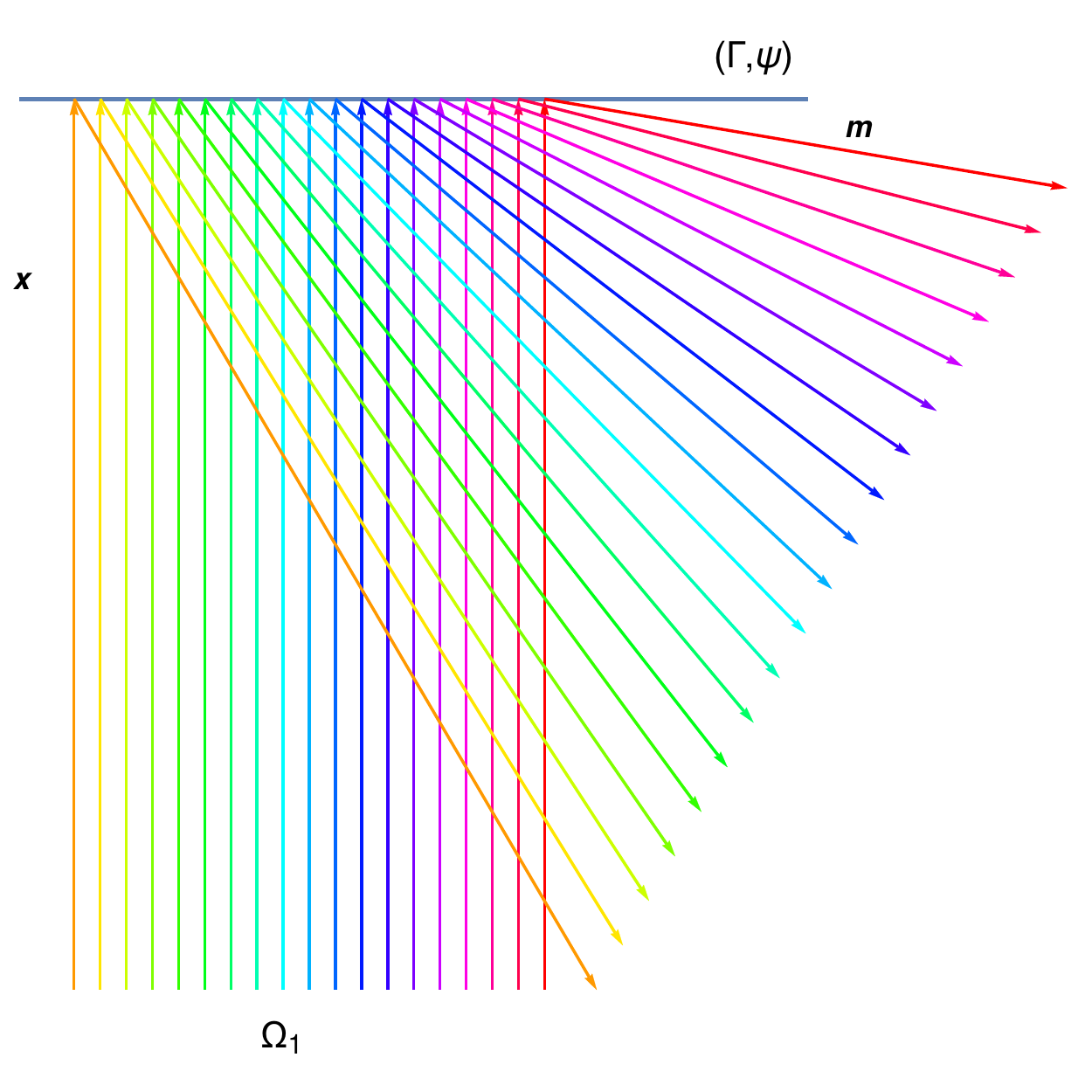}
\caption{Reflection from an extended source (rays are monochromatic; colors are used only for visual purposes).}
\label{fig:general problem}
\end{figure}

\begin{figure}[htbp]
\centering
\includegraphics[width=3in]{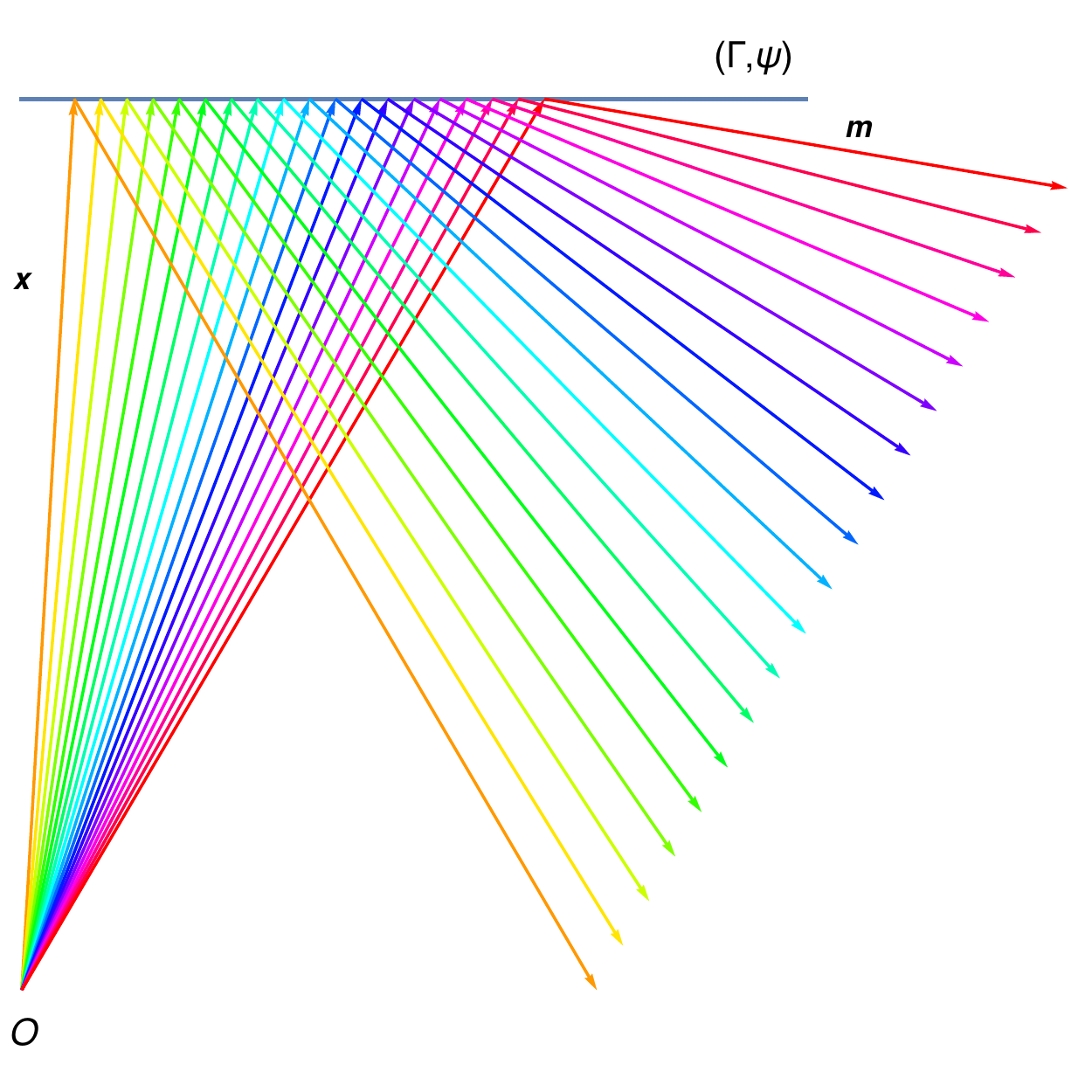}
\caption{Reflection from a point source}
\label{fig:point source reflection}
\end{figure}

In our previous paper \cite{gps} we gave a mathematically rigorous foundation to deal with general metasurfaces and to determine the relationships between the curvature of the surface and the phase discontinuity. 
The present paper builds on that work and shows that phase discontinuities theoretically exist so that there are metalenses that can refract or reflect energy with prescribed energy patterns.
More precisely, a question considered and answered in this paper is the following.
Let $\Omega_1$ be a set on the $x,y$-plane, $\Omega_2$ a set of unit vectors in $\R^3$, and  
$\Gamma=\{z=a\}$ a horizontal plane above the $x,y$-plane. We are given two intensities, i.e, two non negative functions, $f$ defined in $\Omega_1$ and $g$ defined in $\Omega_2$ satisfying the energy conservation condition
\[
\int_{\Omega_1}f(x,y)\,dxdy=\int_{\Omega_2}g(z)\,d\sigma(z),
\] 
where $d\sigma(z)$ denotes as usual the element of area in the unit sphere of $\R^3$.
A collimated beam is emanating from $\Omega_1$. That is, for each $(x,y)\in \Omega_1$ a ray is emitted in the vertical direction $e_3=(0,0,1)$ with intensity $f(x,y)$ and strikes the plane $\Gamma$ at the point $(x,y,a)=P$.
According to the generalized law of reflection \eqref{eq:generalized law of reflection} this ray is reflected by the meta surface $(\Gamma,\psi)$ into a ray having direction $T(x,y)=e_3-\lambda \,e_3-\nabla \psi(x,y,a)$, since $\nu(P)=e_3$. 
The question is then to find a function $\psi$ defined in $\Gamma$, called a phase discontinuity, such that the meta lens, i.e., the pair $(\Gamma,\psi)$
reflects all rays from $\Omega_1$ into $\Omega_2$, that is, $T(\Omega_1)=\Omega_2$, and the energy conservation balance
\begin{equation}\label{eq:conservation of energy over each set}
\int_{E} f(x,y)\,dxdy=\int_{T(E)} g(z)\,d\sigma(z)
\end{equation}
holds for each subset $E$ of $\Omega_1$, see Figure \ref{fig:general problem}.
We show in Section  \ref{subsec:reflection collimated case} that this problem is mathematically solvable.
We also consider a similar question when the incident rays emanate from a point source into a set of unit directions $\Omega_1$, see Figure \ref{fig:point source reflection}. Such a problem is solved in Section \ref{subsec:point source reflection}. 

\begin{figure}[htbp]
\centering
\includegraphics[width=3in]{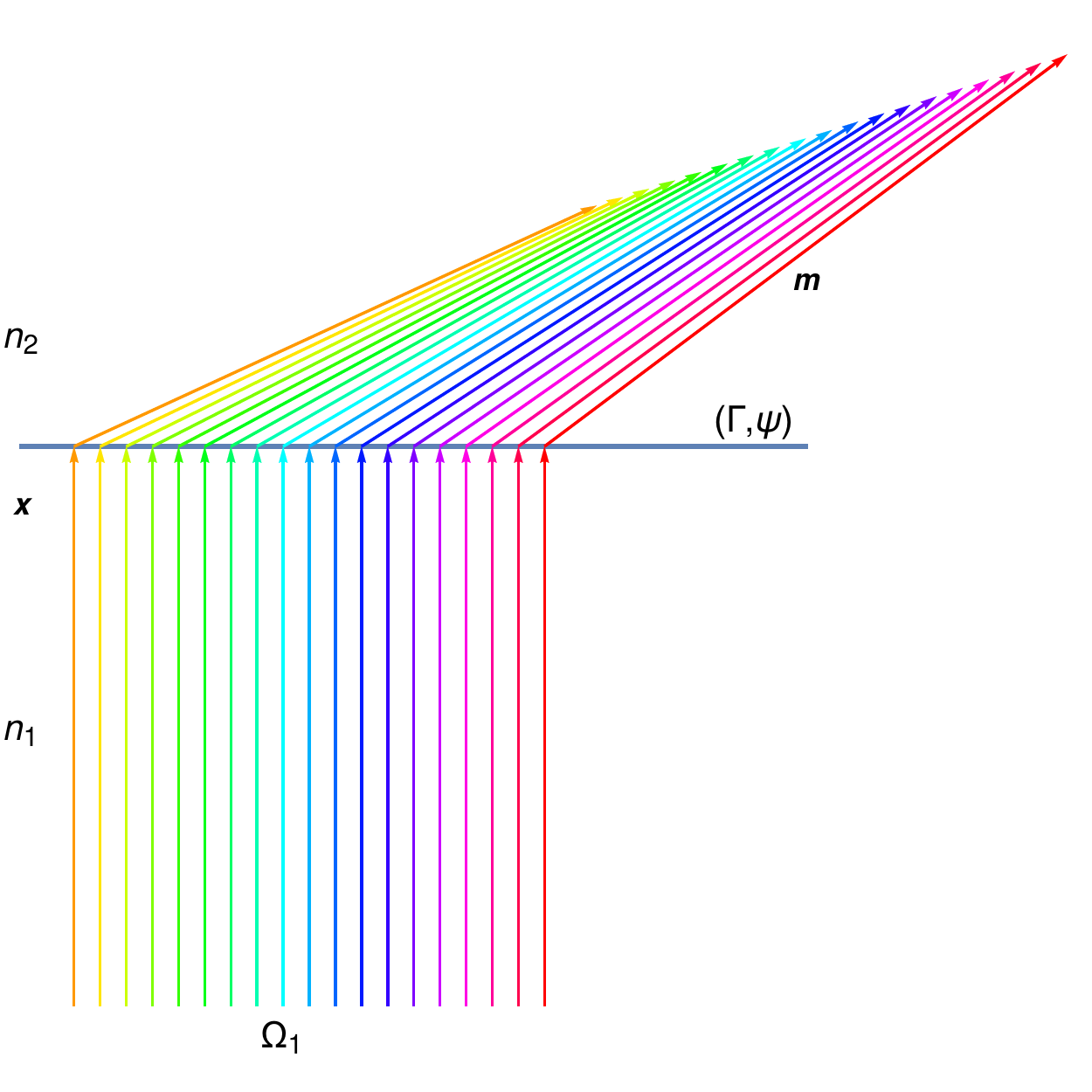}
\caption{Refraction from an extended source}
\label{fig:extended refraction}
\end{figure}

\begin{figure}[htbp]
\centering
\includegraphics[width=3in]{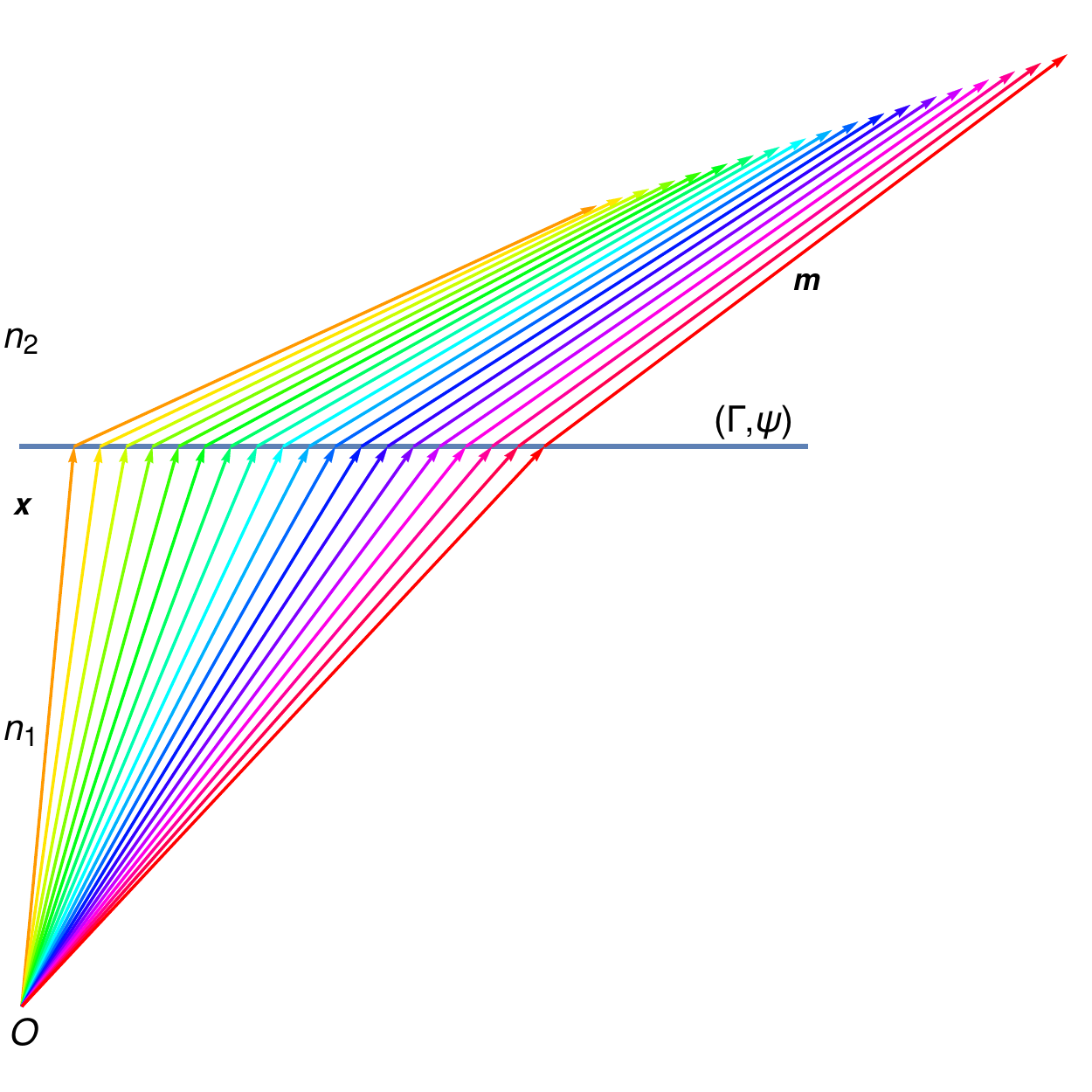}
\caption{Refraction from a point source}
\label{fig:point source refraction}
\end{figure}

In addition, we consider and solve similar problems for refraction using the generalized law \eqref{eq:generalized law of refraction} both in the collimated and point source cases, see Figures \ref{fig:extended refraction} and \ref{fig:point source refraction}, Sections  \ref{subsec:refraction collimated case} and \ref{subsec:refraction point source case}.

It is our purpose to show that each of these problems has a theoretical solution.
To do this we derive the partial differential equation, for each problem, satisfied by the phase discontinuity $\psi$ and show it is a Monge-Amp\`ere equation. Next we show that the resulting equations 
have solutions by application of a result by Urbas  \cite{u}. The equations corresponding to the four problems considered are \eqref{eq:monge ampere reflection}, \eqref{eq:pde for one source far field}, \eqref{eq:eq1} and \eqref{eq:pde for point source refraction}, and they can be regarded as particular cases of \eqref{eq:unique_eq}.

Monge-Amp\`ere equations appear naturally in optics for freeform lens design that have been the subject of recent research, see for example \cite{w}-\cite{gutierrez-sabra:freeformgeneralfields}.
Therefore it is natural that these type of equations appear also for metasurfaces.
Monge-Amp\`ere equations have been recently the subject of important mathematical research due to their connections with various topics such us optimal mass transport. We refer the reader to \cite{Gutierrez-book} and \cite{Figalli-book} for details and references therein.
We hope our contributions in this paper may be useful to understand theoretically what kind of phase discontinuities are possible to design in the applications of metasurfaces.

We mention that recent work using the ideas from \cite{gps} for reflection is done in \cite{biswas-gutierrez-low-2018}
to design graphene-based metasurfaces that can be actively tuned between different regimes of
operation, such as anomalous beam steering and focusing, cloaking, and illusion optics, by applying
electrostatic gating without modifying the geometry of the metasurface.

Finally, if the surface $\Gamma$ is not necessarily a plane, then is possible to derive the corresponding partial differential equation that the phase discontinuity $\psi$ satisfies, in both the reflection and refraction cases. These are equations of Monge-Amp\`ere type that require a more complicated derivation and analysis beyond the scope of this article and
 will appear as a part of the forthcoming PH.D. thesis of the second author.

\subsection{Background}
We begin describing the formulation of the generalized laws of reflection and refraction with phase discontinuities from \cite{gps}.
For refraction consider two homogenous and isotropic media $I$ and $II$ with refractive indices $n_1$ and $n_2$ respectively. Suppose we have an interface surface $\Gamma$ separating media $I$ and $II$, and a function $\psi$, called the phase discontinuity, defined in a small neighborhood of $\Gamma$. If a ray with unit direction $\mathbold x$ emanating from medium $I$ strikes the surface $\Gamma$ at some point $P$, then it is refracted into medium $II$ into a ray having unit direction $\mathbold m$ such that
\begin{equation}\label{eq:generalized law of refraction}
n_1\,\mathbold{x}-n_2\,\mathbold m= \lambda \,\nu(P)+\nabla \psi(P),
\end{equation}  
see \cite[Formula (6)]{gps}, where $\nu(P)$ is the outer unit normal to the surface $\Gamma$ at $P$, $\nabla \psi$ denotes the gradient, and $\lambda$ is a constant depending on $x, \nu(P),\nabla \psi(P)$ and the refractive indices. 
In fact, it is proved in \cite[Formula (11)]{gps} that $\lambda$ can be calculated:
\begin{equation}\label{eq:formula for lambda refraction}
\lambda=\left(n_1\,\mathbold x- \nabla \psi\right)\cdot \nu
-
\sqrt{n_2^2 -|n_1\,\mathbold x-\nabla \psi|^2+ \left(\left(n_1\,\mathbold x-\nabla \psi\right)\cdot \nu \right)^2}.
\end{equation}  
This law is derived in \cite[Section 3]{gps} using wave fronts.
The job of the function $\psi$, concentrated around $\Gamma$, is to change the direction of  the incoming rays. As a difference from the standard Snell law of refraction, the rays here are mainly bent due to the function $\psi$ and not by the change in the refractive indices of the surrounding media. 
In particular, the generalized law of refraction \eqref{eq:generalized law of refraction} makes sense when $n_1=n_2$.
When $\psi$ is constant and $n_1\neq n_2$, we obtain the standard Snell law of refraction.

The case of reflection is when $n_1=n_2$ and since now the reflected vector must be on the same side of $\Gamma$, that  is, $\mathbold{m}\cdot \nu\leq 0$, the generalized reflection law has the form
\begin{align}\label{eq:generalized law of reflection}
\mathbold x-\mathbold m=\lambda/n_1\nu(P)+\nabla (\psi/n_1)(P),
\end{align}
with $\mathbold x$ the unit incident ray, $\mathbold m$ the unit reflected ray, $\nu(P)$ the normal to $\Gamma$ at $P$, and
\begin{equation}\label{eq:formula for lambda reflection}
\lambda=\left(n_1\,\mathbold x- \nabla \psi\right)\cdot \nu
+
\sqrt{n_1^2 -|n_1\,\mathbold x-\nabla \psi|^2+ \left(\left(n_1\,\mathbold x-\nabla \psi\right)\cdot \nu \right)^2},
\end{equation} 
see \cite[Section 3, Remark 1]{gps}.
Once again when $\psi$ is constant this yields the standard reflection law.

\section{Reflection}\label{sec:reflection}
\subsection{Collimated case}\label{subsec:reflection collimated case}
Here we solve the first problem stated in the introduction.
From \eqref{eq:generalized law of reflection} and \eqref{eq:formula for lambda reflection} with $n_1=1$, the vertical ray  emanating from the point $(x,y)\in\Omega_1$ is reflected by the metasurface $(\Gamma,\psi)$ into the unit direction 
\begin{equation}\label{eq:reflection map}
T(x,y)=i(x,y)-\lambda\nu(x,y)-\nabla\psi(x,y),
\end{equation} 
where $i(x,y)=(0,0,1)$,  $\nu(x,y)$ is the normal to $\Gamma=\{z=1\}$, and
\begin{align*}
\lambda=(i-\nabla\psi)\cdot \nu+ \sqrt{1-\big(|i-\nabla \psi|^2-[(i-\nabla\psi)\cdot \nu]^2\big)}
=1+\sqrt{1-\psi^2_x-\psi^2_y}.
\end{align*} 
We remark that  in the last identity we have used that $\psi_z=0$ since we seek a phase discontinuity $\psi$ tangential to the surface $\Gamma$. 

Therefore 
\begin{align}\label{eq:formula for reflection T(x,y)}
T(x,y)=(T_1,T_2,T_3)
=-\left(\psi_x(x,y),\psi_y(x,y), \sqrt{1-\psi_x^2(x,y)-\psi_y^2(x,y)}\right)\nonumber.
\end{align}
From the conservation of energy condition \eqref{eq:conservation of energy over each set} and the formula of change of variables for surface integrals
\begin{equation}\label{eq:definition of solution aleksandrov}
\int_Ef(x)\,dx=\int_{T(E)}g(y)\,d\sigma(y)=\int_Eg(T(z))|J_T(z)|\,dz,
\end{equation}
for each open set $E\subset \Omega_1$, and 
where $|J_T|=|T_x(x,y)\times T_y(x,y)|$. 
From \eqref{eq:definition of solution aleksandrov} we obtain
\begin{equation}\label{eq:differential equation with Jacobian collimated}
f(x)=g(T(x,y))\,|J_T(x,y)|\qquad \text{for $(x,y)\in \Omega_1$}.
\end{equation}
To calculate $|J_T(x,y)|$, since $|T(x,y)|=1$, differentiating with respect to $x$ and $y$ yields the equations 
$T\cdot T_x=T\cdot T_y=0$.
Hence, assuming $T_3(x,y)\neq 0$ and solving these equations in $(T_3)_x$ and $(T_3)_y$ we get 
\begin{align*}
(T_3)_x = - \frac{T_1(T_1)_x+T_2(T_2)_x}{T_3} \quad \text{and} \quad (T_3)_y = - \frac{T_1(T_1)_y+T_2(T_2)_y}{T_3}. 
\end{align*}
From an elementary calculation with the determinant we obtain
\begin{equation*}
T_x\times T_y=\frac{1}{T_3}\det\left(
\begin{matrix}
(T_1)_x& (T_1)_y\\
(T_2)_x& (T_2)_y
\end{matrix}\right)\,T.
\end{equation*}
Hence
\begin{equation*}
|J_T|=\frac{1}{|T_3(x,y)|}|\det(D^2\psi)|,
\end{equation*}
where $D^2\psi$ is the matrix of the second derivatives in $x$ and $y$.
Therefore from \eqref{eq:differential equation with Jacobian collimated} the phase discontinuity $\psi$ satisfies the following Monge-Amp\`ere equation
\begin{equation}\label{eq:monge ampere reflection}
\frac{1}{\sqrt{1-\psi^2_x(x,y)-\psi^2_y(x,y)}}|\det(D^2\psi)|=\frac{f(x,y)}{g\left(T(x,y)\right)}.
\end{equation}

To show that \eqref{eq:monge ampere reflection} has solutions we invoke \cite[Theorem 2]{u}, which in simpler terms says the following: {\it if $D_1,D_2$ are uniformly convex smooth domains in $\R^n$, $f_1>0$ is a smooth function in $D_1$, $f_2>0$ is a smooth function in $D_2$, satisfying 
\begin{equation}\label{hyp1}
\int_{D_1} f_1(x)\,dx=\int_{D_2}f_2(p)\,dp,
\end{equation}
then the boundary value problem
\begin{equation*}
\det(D^2u)=\frac{f_1(x)}{f_2(\nabla u)} \quad \text{in} \quad D_1, \quad \nabla u(D_1)=D_2,
\end{equation*}
has a convex $C^2$ solution, and any two such solutions differ by a constant.}

In fact, to apply this result to show existence of solutions to \eqref{eq:monge ampere reflection}, set $n=2$, let
\begin{align*}
&f_1(x,y)=f(x,y)\text{ for $(x,y)\in D_1=\Omega_1$ }, \\
&f_2(p_1,p_2)= \frac{g\left(-\left(p_1,p_2,\sqrt{1-p_1^2-p_2^2}\right)\right)}{\sqrt{1-p_1^2-p_2^2}},
\end{align*}
for  $(p_1,p_2)\in D_2=-\Pi(\Omega_2)$, where $\Pi$ is the orthogonal projection from a set on the unit sphere onto the $x,y$-plane. In particular, $\Omega_2$ is a subset of the lower unit hemisphere $z\leq 0$.
We need to verify \eqref{hyp1}.
From the conservation of energy assumption
\begin{equation*}
\int_{\Omega_1}f(x)\, dx=\int_{\Omega_2}g(y)\, d\sigma(y),
\end{equation*}
and using the parametrization $q=(q_1,q_2)\to\left(q,-\sqrt{1-|q|^2}\right)$ we can write
\begin{align*}
\int_{\Omega_2}g(y)\, d\sigma(y)&=\int_{\Pi(\Omega_2)}\frac{g\left(q,-\sqrt{1-|q|^2}\right)}{\sqrt{1-|q|^2}}\,dq\\
&= 
\int_{-\Pi(\Omega_2)}\frac{g\left(-p,-\sqrt{1-|p|^2}\right)}{\sqrt{1-|p|^2}}\,dp 
= 
\int_{D_2}f_2(p_1,p_2)\,dp.
\end{align*}
Therefore \eqref{hyp1} holds and hence the existence of solutions to  \eqref{eq:monge ampere reflection} follows.

\subsection{Point Source Reflection}\label{subsec:point source reflection}

We now have a domain $\Omega_1$ of the unit sphere in $\R^3$ and rays emanate from the origin with intensity $f(x)\geq 0$ for each $x\in \Omega_1$.
$\Omega_2$ is as in the previous section, i.e., a domain of the unit sphere and  
 $g> 0$ is a function in $\Omega_2$ such that the following energy conservation condition holds
\begin{equation}\label{eq:conservation of energy whole domains}
\int_{\Omega_1}f(x)\, d\sigma(x)=\int_{\Omega_2}g(y)\, d\sigma(y).
\end{equation}
Again $\Gamma$ is the plane $z=1$. Of course, we assume that rays from the origin with unit direction in $\Omega_1$ reach the plane $\Gamma$.
The question is then to find a phase discontinuity $\psi$ on $\Gamma$ such that all rays emitted from the origin with direction $x\in \Omega_1$ and intensity $f(x)$, are reflected by the meta surface $(\Gamma,\psi)$ into $\Omega_2$ such that 
\begin{equation}\label{consen}
\int_{E} f(x)\,d\sigma(x)=\int_{T(E)} g(y)\,d\sigma(y)
\end{equation}
for each subset $E$ of $\Omega_1$ and $T(\Omega_1)=\Omega_2$, where $T$ is the reflection map.
In order to find the equation $\psi$ satisfies, we parametrize the domains in the sphere using spherical coordinates: $s(u,v)=(\cos{u}\sin{v}, \sin{u}\sin{v}, \cos{v})$, $0\leq u\leq 2\pi$, $0\leq v\leq \pi/2$. Parametrizing $\Omega_1$ in these coordinates we obtain
$\Omega_1=s(O)$, for some domain $O \subset [0,2\pi]\times [0,\pi/2]$. Re writing the integrals in \eqref{consen} in spherical coordinates, and letting $s(U)=E$, we have
\begin{align*}
&\int_U f(s(u,v))|s_u\times s_v| dudv=\int_{E} f(x)\,d\sigma(x)\\
&=\int_{T(E)} g(y)\,d\sigma(y)=\int_U g(T(s(u,v)))|(T\circ s)_u\times(T\circ s)_v|dudv.
\end{align*} 
Since this equation must hold for all open sets $U\subset O$, it follows that $T$ satisfies the equation
\begin{equation}\label{diffeq}
\frac{|(T\circ s)_u\times(T\circ s)_v|}{|s_u\times s_v|}=\frac{f(s(u,v))}{g(T(s(u,v)))}.
\end{equation}

The plane $\Gamma$ is described in spherical coordinates by the polar radius
\begin{equation}\label{eq:polar radius for plane z=1}
r(u,v)=\frac{1}{\cos{v}}s(u,v)=(\cos{u}\tan{v}, \sin{u}\tan{v}, 1).
\end{equation}
From \eqref{eq:generalized law of reflection} with $n_1=1$, if the incident ray has direction $i=s(u,v)$, then the reflected ray that has unit direction
\begin{align*}
T(s(u,v))&=s(u,v)-\lambda\nu-\nabla \psi(r(u,v)),
\end{align*} 
where $\nu=(0,0,1)$ is the normal to $\Gamma$ at the incident point. Since we seek, as before, for a phase $\psi$ tangential to $\Gamma$, we have $\nabla\psi(x,y,1)=(\psi_x(x,y,1),\psi_y(x,y,1),0)$. 
In addition, from \eqref{eq:formula for lambda reflection}
\begin{align*}
\lambda &=i \cdot \nu + \sqrt{1-(|i-\nabla \psi|^2-(i\cdot \nu)^2)}\\
&=\cos{v} + \sqrt{1-(\cos{u}\sin{v}-\psi_x(r(u,v)))^2 -(\sin{u}\sin{v}-\psi_y(r(u,v)))^2}\\
&=\cos{v}+\sqrt{\Delta},
\end{align*}
where $\Delta=1-(\cos{u}\sin{v}-\psi_x(r(u,v)))^2-(\sin{u}\sin{v}-\psi_y(r(u,v)))^2$.
Therefore writing $T$ in components
\begin{align}\label{eq:components of T reflection one source}
&T(s(u,v))=\left(T_1(s(u,v)),T_2(s(u,v)),T_3(s(u,v))\right)\notag \\
&=\left(\cos{u}\sin{v}-\psi_x(r(u,v)), \sin{u}\sin{v}-\psi_y(r(u,v)),-\sqrt{\Delta}\right).
\end{align}
Since $|T(s(u,v))|=1$, it follows as in Section \ref{subsec:reflection collimated case} that 
\begin{equation}\label{eq:formula cross product T spherical}
|(T\circ s)_u\times (T\circ s)_v|=\frac{1}{|T_3\circ s|}\vline\det
\begin{pmatrix}
(T_1\circ s)_u & (T_1\circ s)_v \\
(T_2\circ s)_u & (T_2\circ s)_v
\end{pmatrix}\vline.
\end{equation}
From \eqref{eq:components of T reflection one source}
\begin{align*}
(T_1\circ s)_u&=-\sin{u}\sin{v}-\psi_{xx}(r(u,v))(-\sin{u}\tan{v})-\psi_{xy}(r(u,v))(\cos{u}\tan{v}),\\
(T_1\circ s)_v&=\cos{u}\cos{v}-\psi_{xx}(r(u,v))\left(\frac{\cos{u}}{\cos^2{v}}\right)-\psi_{xy}(r(u,v))\left(\frac{\sin{u}}{\cos^2{v}}\right),\\
(T_2\circ s)_u&=\cos{u}\sin{v}-\psi_{xy}(r(u,v))(-\sin{u}\tan{v})-\psi_{yy}(r(u,v))(\cos{u}\tan{v}),\\
(T_2\circ s)_v&=\sin{u}\cos{v}-\psi_{xy}(r(u,v))\left(\frac{\cos{u}}{\cos^2{v}}\right)-\psi_{yy}(r(u,v))\left(\frac{\sin{u}}{\cos^2{v}}\right).
\end{align*}
Inserting these in \eqref{eq:formula cross product T spherical} yields
\begin{align}\label{eq:formula absolute cross product spherical}
|(T\circ s)_u\times (T\circ s)_v|=\frac{1}{|T_3\circ s|}\vline\det\big(A(u,v) - D^2_{(x,y)}\psi(r(u,v)) B(u,v) \big)\vline 
\end{align}
where 
\begin{align*}
A(u,v)&=
\begin{pmatrix}
-\sin{u}\sin{v} & \cos{u}\cos{v}\\
\cos{u}\sin{v} & \sin{u}\cos{v}
\end{pmatrix},
\\
B(u,v)&=
\begin{pmatrix}
-\sin{u}\tan{v} & \frac{\cos{u}}{\cos^2{v}}\\
\cos{u}\tan{v} & \frac{\sin{u}}{\cos^2{v}}
\end{pmatrix}.
\end{align*}
We can re write the above quantities in rectangular coordinates
noticing that $x=\cos{u}\tan v$, $y=\sin{u}\tan v$, $r(u,v)=(x,y,1)$, $\sqrt{x^2+y^2+1}=\dfrac{1}{\cos{v}}$ and $\sqrt{x^2+y^2}=\tan{v}$.
We obtain 
\begin{align*}
T_1&=\frac{x}{\sqrt{x^2+y^2+1}}-\psi_x(x,y,1)=\left(\sqrt{x^2+y^2+1}-\psi(x,y,1)\right)_x,\\
T_2&=\frac{y}{\sqrt{x^2+y^2+1}}-\psi_y(x,y,1)=\left(\sqrt{x^2+y^2+1}-\psi(x,y,1)\right)_y,\\
|T_3|&=\sqrt{1-\left(\frac{x}{\sqrt{x^2+y^2+1}}-\psi_x(x,y,1)\right)^2-\left(\frac{y}{\sqrt{x^2+y^2+1}}-\psi_y(x,y,1)\right)^2}\\   
&=\sqrt{1-\left(\left(\sqrt{x^2+y^2+1}-\psi(x,y,1)\right)_x\right)^2-\left( \left(\sqrt{x^2+y^2+1}-\psi(x,y,1)\right)_y\right)^2},\\
A&=
\begin{pmatrix}
\dfrac{-y}{\sqrt{x^2+y^2+1}} & \dfrac{x}{\sqrt{x^2+y^2+1}\sqrt{x^2+y^2}}\\
\dfrac{x}{\sqrt{x^2+y^2+1}} & \dfrac{y}{\sqrt{x^2+y^2+1}\sqrt{x^2+y^2}}
\end{pmatrix},\\
B&=
\begin{pmatrix}
-y & \dfrac{x(1+x^2+y^2)}{\sqrt{x^2+y^2}}\\
x & \dfrac{y(1+x^2+y^2)}{\sqrt{x^2+y^2}}
\end{pmatrix}.
\end{align*}
Also,
\begin{align}\label{eq:su cross sv}
|s_u\times s_v|&=\frac{1}{|\cos v|}\vline-\sin^2 u\cos{v}\sin{v}-\cos^2 u\cos{v}\sin{v}\vline\\
&=\sin{v}=\frac{\sqrt{x^2+y^2}}{\sqrt{x^2+y^2+1}}.\nonumber
\end{align}
Now notice that 
\begin{align*}
\det(A-D^2\psi B)=\det(B)\det(AB^{-1}- D^2\psi),
\end{align*}
with
\begin{align}\label{eq:matrix B^-1 reflection}
B^{-1}=-\frac{1}{(x^2+y^2+1)\sqrt{x^2+y^2}}\begin{pmatrix}
\dfrac{y(1+x^2+y^2)}{\sqrt{x^2+y^2}} & \dfrac{-x(1+x^2+y^2)}{\sqrt{x^2+y^2}}\\
-x & -y
\end{pmatrix}
\end{align}
and
\begin{align*}
&AB^{-1}=\begin{pmatrix}
\dfrac{y^2}{b(x)}+\dfrac{x^2}{c(x)} & 
\dfrac{-xy}{b(x)}+\dfrac{xy}{c(x)}\\
\dfrac{-xy}{b(x)}+\dfrac{xy}{c(x)} & 
\dfrac{x^2}{b(x)}+\dfrac{y^2}{c(x)} 
\end{pmatrix}=D^2\Big(\sqrt{x^2+y^2+1}\Big),
\end{align*}
where $b(x)=(x^2+y^2)^{}(x^2+y^2+1)^{1/2}$ and $c(x)=(x^2+y^2)^{}(x^2+y^2+1)^{3/2}$.
Therefore
\begin{align*}
\det(A-D^2\psi B)=\det(B)\det\left(D^2\left(\sqrt{x^2+y^2+1}-\psi\right)\right).
\end{align*}
Letting $\phi(x,y)=\sqrt{x^2+y^2+1}-\psi(x,y)$, using the last equation in \eqref{eq:formula absolute cross product spherical}, and using \eqref{eq:su cross sv}, we obtain from \eqref{diffeq} that $\phi$ satisfies the following equation
\begin{align}\label{eq:pde for one source far field}
&\frac{(x^2+y^2+1)^{3/2}}{\sqrt{1-\phi_x^2(x,y)-\phi_y^2(x,y)}}|\det(D^2\phi(x,y))|\\
&=
\frac{f\left(\frac{1}{\sqrt{x^2+y^2+1}}(x,y,1)\right)}{g\left(\phi_x(x,y),\phi_y(x,y),-\sqrt{1-\phi_x^2(x,y)-\phi_y^2(x,y)}\right)}. \nonumber
\end{align}
The above equation holds for $\left(x,y\right)\in D$ where $D$ is obtained as follows: for each direction $e\in \Omega_1$, the ray with this direction intersect the plane $z=1$ at a unique point $(x,y)$, this collection of $x$ and $y$ is $D$.  

We now proceed as in the previous section to show existence of solutions to \eqref{eq:pde for one source far field}.
To this end we need to identity the functions $f_1,f_2$ in \eqref{hyp1}.
Parametrizing $\Omega_1$ by $q:D\to \Omega_1$ with $q(x,y)=\frac{1}{\sqrt{x^2+y^2+1}}(x,y,1)$, we let 
$f_1(x,y)=\dfrac{f\left( q(x,y)\right)}{(x^2+y^2+1)^{3/2}}$ for $(x,y)\in D_1=D$. 
Also let $f_2(p_1,p_2)=\dfrac{g\left(p_1,p_2,-\sqrt{1-p_1^2-p_2^2} \right)}{\sqrt{1-p_1^2-p_2^2}}$, for $(p_1,p_2)\in D_2=\Pi(\Omega_2)$. With these choices and observing that 
\begin{align*}
\int_{\Omega_1}f(z) \, d\sigma(z) 
= 
\int_D \frac{f\left(\frac{1}{\sqrt{x^2+y^2+1}}(x,y,1)\right)}{(x^2+y^2+1)^{3/2}}\, dxdy,
\end{align*}
a similar calculation as at the end of last section, we obtain that \eqref{eq:conservation of energy whole domains} is equivalent to \eqref{hyp1} and therefore existence of solutions to \eqref{eq:pde for one source far field} follows as before invoking \cite[Theorem 2]{u}.

\section{Refraction}\label{sec:refraction two cases}
Here we solve two problems similar to the ones considered in the previous section but for refraction.
\subsection{Collimated case}\label{subsec:refraction collimated case}
Incident rays are emitted from an open set $\Omega_1$ of the $x$-$y$ plane with direction $i(x,y)=e_3=(0,0,1)$, and $\Gamma$ is the plane $z=1$.

From the generalized law of refraction \eqref{eq:generalized law of refraction} and \eqref{eq:formula for lambda refraction}, the meta surface $(\Gamma,\psi)$ refracts the incident ray $i(x,y)$ into a ray $r(x,y)$ with direction satisfying
$$
n_1i(x,y)-n_2r(x,y)=\lambda\nu(x,y) +\nabla\psi(x,y),
$$
where $n_1$ and $n_2$ are the refractive indices of the two homogeneous and isotropic media separated by the plane $\Gamma$, $\nu$ is the unit normal to the plane $\Gamma$. Also 
\begin{align*}
\lambda &= (n_1i-\nabla\psi)\cdot \nu-\sqrt{n_2^2-|n_1i-\nabla\psi|^2+[(n_1i-\nabla\psi)\cdot\nu]^2}\\
&=
n_1-\sqrt{n_2^2-(\psi_x^2+\psi_y^2)},
\end{align*}
since we seek $\psi$ tangential to $\Gamma$; i.e., $\psi_z=0$.
We then let $T:\Omega_1\to\Omega_2$ to be
\begin{align*}
T(x,y):=r(x,y)=\left(-\frac{1}{n_2}\psi_x(x,y),-\frac{1}{n_2}\psi_y(x,y),\sqrt{1-\frac{1}{n_2^2}\left(\psi_x^2(x,y)+\psi_y^2(x,y)\right)}\right).  
\end{align*}
We seek $\psi$ defined on $\Gamma$ with $T(\Omega_1)=\Omega_2$ and 
satisfying the conservation of energy balance 
\begin{align*}
\int_Ef(x)\,dx=\int_{T(E)}g(y)\,d\sigma(y)=\int_Eg(T(z))|J_T|\,dz\quad \text{for each}\quad E\subset \Omega_1,
\end{align*}
where $|J_T|=|T_x(x,y)\times T_y(x,y)|$. Since $|T(x,y)|=1$ and similarly as for reflection, we have that
\begin{equation*}
|J_T|=\frac{1}{|T_3(x,y)|}\left|\det\left(D^2\frac{1}{n_2}\psi\right)\right|.
\end{equation*}
Therefore proceeding as in the reflection case, the phase discontinuity $\psi$ must satisfy the following Monge-Amp\'ere equation
\begin{align}
\label{eq:eq1}
\frac{1}{\sqrt{1-\frac{1}{n_2^2}\left(\psi_x^2(x,y)+\psi_y^2(x,y)\right)}}\left|\det\left(D^2\frac{1}{n_2}\psi\right)\right|=\frac{f(x,y)}{g\left(-\frac{1}{n_2}\psi_x,-\frac{1}{n_2}\psi_y,\sqrt{1-\frac{1}{n_2^2}\left(\psi_x^2+\psi_y^2\right)}\right)}; 
\end{align}
notice that this equation is independent of the value of $n_1$.
Similar to the reflection case, $T(\Omega_1)=\Omega_2$ implies that $\frac{1}{n_2}(\psi_x,\psi_y)\in -\Pi(\Omega_2)$ where $\Pi$ is once again the orthogonal projection onto the $x$-$y$ plane.
We claim, also in this case, that \cite[Theorem 2]{u} can be applied to obtain a solutions $\psi$ to \eqref{eq:eq1}.
Indeed, letting 
\begin{align*}
&f_1(x,y)=f(x,y)\text{ for $(x,y)\in D_1=\Omega_1$ },\\
&f_2(p_1,p_2)= \frac{g\left(-\frac{1}{n_2}p_1,-\frac{1}{n_2}p_2,\sqrt{1-\frac{1}{n_2^2}\left(p_1^2+p_2^2\right)}\right)}{\sqrt{1-\frac{1}{n_2^2}\left(p_1^2+p_2^2\right)}},
\end{align*}
for $(p_1,p_2)\in D_2=-n_2\Pi(\Omega_2)$,
and proceeding as before we obtain that 
\begin{equation*}
\int_{\Omega_1}f(x)\, dx=\int_{\Omega_2}g(y)\, d\sigma(y).
\end{equation*}
is equivalent to \eqref{hyp1} and so existence of solutions follows as before.

\subsection{Point Source Refraction}\label{subsec:refraction point source case}

We now analyze a problem similar to the one in Section  \ref{subsec:point source reflection} for refraction.
That is, rays emanate for a point source and we seek a phase discontinuity $\psi$ defined on the plane $\Gamma=\{z=1\}$ so that the refraction map $T$ (to be calculated in a moment) satisfies  
the conservation of energy condition (\ref{consen}). As in Section  \ref{subsec:point source reflection}, this implies (\ref{diffeq}), i.e., 
\[
\frac{|(T\circ s)_u\times(T\circ s)_v|}{|s_u\times s_v|}=\frac{f(s(u,v))}{g(T(s(u,v)))},
\]
and $T(\Omega_1)=\Omega_2$.
Let us calculate the refraction map $T$.
As in Section \ref{subsec:point source reflection}, the plane $\Gamma$ is described by the polar radius 
\eqref{eq:polar radius for plane z=1}.
Then from \eqref{eq:generalized law of refraction}, the refracted ray has unit direction
\begin{align*}
T(s(u,v))&=\frac{n_1}{n_2}s(u,v)-\frac{1}{n_2}\lambda\nu-\frac{1}{n_2}\nabla \psi(r(u,v)),
\end{align*} 
where $\nu=(0,0,1)$ is the normal to $\Gamma$ at the incident point, $s(u,v)$ are spherical coordinates, $\nabla\psi(x,y,1)=(\psi_x(x,y,1),\psi_y(x,y,1),0)$ (since we seek a phase discontinuity $\psi$ tangential to $\Gamma$), and 
\begin{align*}
\lambda &=n_1i \cdot \nu + \sqrt{n_2^2-(|n_1i-\nabla \psi|^2-(n_1i\cdot \nu)^2)}\\
&=n_1\cos{v} + \sqrt{n_2^2-(n_1\cos{u}\sin{v}-\psi_x(r(u,v)))^2-(n_1\sin{u}\sin{v}-\psi_y(r(u,v)))^2}\\
&=n_1\cos{v}+\sqrt{\Delta},
\end{align*}
where $\Delta=n_2^2-(n_1\cos{u}\sin{v}-\psi_x(r(u,v)))^2-(n_1\sin{u}\sin{v}-\psi_y(r(u,v)))^2$.
Therefore 
\begin{align*}
T(s(u,v))&=\left(T_1(s(u,v)),T_2(s(u,v)),T_3(s(u,v))\right) \\
&=\(\frac{n_1}{n_2}\cos{u}\sin{v}-\frac{1}{n_2}\psi_x(r(u,v)), \frac{n_1}{n_2}\sin{u}\sin{v}-\frac{1}{n_2}\psi_y(r(u,v)),\frac{1}{n_2}\sqrt{\Delta}\).
\end{align*}
Since $|T(s(u,v))|=1$ we have as in \eqref{eq:formula cross product T spherical} that 
\begin{equation}\label{eq:formula cross product T spherical refraction}
|(T\circ s)_u\times (T\circ s)_v|=\frac{1}{|T_3\circ s|}\left|\det
\begin{pmatrix}
(T_1\circ s)_u & (T_1\circ s)_v \\
(T_2\circ s)_u & (T_2\circ s)_v
\end{pmatrix}\right|.
\end{equation}
On the other hand,
\begin{align*}
(T_1\circ s)_u&=-\frac{n_1}{n_2}\sin{u}\sin{v}-\frac{1}{n_2}\psi_{xx}(r(u,v))(-\sin{u}\tan{v})-\frac{1}{n_2}\psi_{xy}(r(u,v))(\cos{u}\tan{v}),\\
(T_1\circ s)_v&=\frac{n_1}{n_2}\cos{u}\cos{v}-\frac{1}{n_2}\psi_{xx}(r(u,v))\left(\frac{\cos{u}}{\cos^2{v}}\right)-\frac{1}{n_2}\psi_{xy}(r(u,v))\left(\frac{\sin{u}}{\cos^2{v}}\right),\\
(T_2\circ s)_u&=\frac{n_1}{n_2}\cos{u}\sin{v}-\frac{1}{n_2}\psi_{xy}(r(u,v))(-\sin{u}\tan{v})-\frac{1}{n_2}\psi_{yy}(r(u,v))(\cos{u}\tan{v}),\\
(T_2\circ s)_v&=\frac{n_1}{n_2}\sin{u}\cos{v}-\frac{1}{n_2}\psi_{xy}\left(r(u,v))(\frac{\cos{u}}{\cos^2{v}}\right)-\frac{1}{n_2}\psi_{yy}(r(u,v))\left(\frac{\sin{u}}{\cos^2{v}}\right).
\end{align*}
Inserting these into \eqref{eq:formula cross product T spherical refraction} yields
\begin{align*}
|(T\circ s)_u\times (T\circ s)_v|=\frac{1}{|T_3\circ s|}\vline\det\big(A(u,v) -\frac{1}{n_2} D^2_{(x,y)}\psi(r(u,v)) B(u,v) \big)\vline
\end{align*}
where 
\begin{align*}
&A(u,v)=\frac{n_1}{n_2}
\begin{pmatrix}
-\sin{u}\sin{v} & \cos{u}\cos{v}\\
\cos{u}\sin{v} & \sin{u}\cos{v}
\end{pmatrix},
\\
&B(u,v)=
\begin{pmatrix}
-\sin{u}\tan{v} & \frac{\cos{u}}{\cos^2{v}}\\
\cos{u}\tan{v} & \frac{\sin{u}}{\cos^2{v}}
\end{pmatrix}.
\end{align*}
As in the point source reflection case, Section \ref{subsec:point source reflection}, we can re write the above quantities in rectangular coordinates
noticing that $x=\cos{u}\tan v$, $y=\sin{u}\tan v$, $r(u,v)=(x,y,1)$, $\sqrt{x^2+y^2+1}=\dfrac{1}{\cos{v}}$ and $\sqrt{x^2+y^2}=\tan{v}$.
We obtain 
\begin{align*}
T_1=&\frac{n_1}{n_2}\frac{x}{\sqrt{x^2+y^2+1}}-\frac{1}{n_2}\psi_x(x,y,1)
=\left(\frac{n_1}{n_2}\sqrt{x^2+y^2+1}-\frac{1}{n_2}\psi(x,y,1)\right)_x,\\
T_2=&\frac{n_1}{n_2}\frac{y}{\sqrt{x^2+y^2+1}}-\frac{1}{n_2}\psi_y(x,y,1)
=\left(\frac{n_1}{n_2}\sqrt{x^2+y^2+1}-\frac{1}{n_2}\psi(x,y,1)\right)_y,\\
|T_3|=&\sqrt{1-\left(\frac{n_1}{n_2}\frac{x}{\sqrt{x^2+y^2+1}}-\frac{1}{n_2}\psi_x(x,y,1)\right)^2-\left(\frac{n_1}{n_2}\frac{y}{\sqrt{x^2+y^2+1}}-\frac{1}{n_2}\psi_y(x,y,1)\right)^2}\\   
=&\sqrt{1-\left(\left(\frac{n_1}{n_2}\sqrt{x^2+y^2+1}-\frac{1}{n_2}\psi(x,y,1)\right)_x\right)^2-\left( \left(\frac{n_1}{n_2}\sqrt{x^2+y^2+1}-\frac{1}{n_2}\psi(x,y,1)\right)_y\right)^2},\\
A=&\frac{n_1}{n_2}
\begin{pmatrix}
\dfrac{-y}{\sqrt{x^2+y^2+1}} & \dfrac{x}{\sqrt{x^2+y^2+1}\sqrt{x^2+y^2}}\\
\dfrac{x}{\sqrt{x^2+y^2+1}} & \dfrac{y}{\sqrt{x^2+y^2+1}\sqrt{x^2+y^2}}
\end{pmatrix},\\
B=&
\begin{pmatrix}
-y & \dfrac{x(1+x^2+y^2)}{\sqrt{x^2+y^2}}\\
x & \dfrac{y(1+x^2+y^2)}{\sqrt{x^2+y^2}}
\end{pmatrix}.
\end{align*}
Also from \eqref{eq:su cross sv}
$
|s_u\times s_v|
=\frac{\sqrt{x^2+y^2}}{\sqrt{x^2+y^2+1}}.
$
Now notice that 
\begin{align*}
\det\left(A-\frac{1}{n_2}D^2\psi B\right)=\det(B)\det\left(AB^{-1}- \frac{1}{n_2}D^2\psi\right),
\end{align*}
with $B^{-1}$ as in \eqref{eq:matrix B^-1 reflection}

and
\begin{align*}
&AB^{-1}=\frac{n_1}{n_2}\footnotesize{\begin{pmatrix}
\dfrac{y^2}{b(x)}+\dfrac{x^2}{c(x)} & 
\dfrac{-xy}{b(x)}+\dfrac{xy}{c(x)}\\
\dfrac{-xy}{b(x)}+\dfrac{xy}{c(x)} & 
\dfrac{x^2}{b(x)}+\dfrac{y^2}{c(x)} 
\end{pmatrix}}=\frac{n_1}{n_2}D^2\Big(\sqrt{x^2+y^2+1}\Big),
\end{align*}
where $b(x)=(x^2+y^2)^{}(x^2+y^2+1)^{1/2}$ and $c(x)=(x^2+y^2)^{}(x^2+y^2+1)^{3/2}$.
Therefore
\begin{align*}
\det\left(A-\frac{1}{n_2}D^2\psi  B\right)=\det(B)\det\left(D^2\left(\frac{n_1}{n_2}\sqrt{x^2+y^2+1}-\frac{1}{n_2}\psi\right)\right).
\end{align*}
Letting $\phi(x,y)=\frac{n_1}{n_2}\sqrt{x^2+y^2+1}-\frac{1}{n_2}\psi(x,y)$, we obtain that $\phi$ satisfies the following equation
\begin{align}\label{eq:pde for point source refraction}
&\frac{(x^2+y^2+1)^{3/2}}{\sqrt{1-\phi_x^2(x,y)-\phi_y^2(x,y)}}|\det(D^2\phi(x,y))|
\\&=
\frac{f\left(\frac{1}{\sqrt{x^2+y^2+1}}(x,y,1)\right)}{g\left(\phi_x(x,y),\phi_y(x,y),\sqrt{1-\phi_x^2(x,y)-\phi_y^2(x,y)}\right)}\nonumber.
\end{align}
The above equation holds for $\left(x,y\right)\in D$ where $D$ is obtained as at the end of Section \ref{subsec:point source reflection}.
Existence of solutions to this equation follows as before letting 
$f_1(x,y)=\dfrac{f\left( \frac{1}{\sqrt{x^2+y^2+1}}(x,y,1)\right)}{(x^2+y^2+1)^{3/2}}$ for $(x,y)\in D_1=D$,
and $f_2(p_1,p_2)=\dfrac{g\left(p_1,p_2,\sqrt{1-p_1^2-p_2^2} \right)}{\sqrt{1-p_1^2-p_2^2}}$ for $(p_1,p_2)\in D_2=\Pi(\Omega_2)$, where $\Pi$ is once again the orthogonal projection.

\begin{remark}\rm
If a ray is emitted from a point $Q$ and strikes the plane $\Gamma=\{z=1\}$ at the point $P=(x,y,1)$, let $d_Q(x,y)$ be the distance from $Q$ to $P$. 
In the collimated case, since all rays are vertical $d_Q(x,y)=1$.
And when the point source $Q$ is the origin, $d_Q(x,y)=\sqrt{x^2+y^2+1}$.
Then writing  $\phi(x,y)=\dfrac{n_1}{n_2}d_{Q}(x,y)-\dfrac{1}{n_2}\psi(x,y)$, and noticing that $n_1=n_2=1$ in the reflection cases, 
%
the equations \eqref{eq:monge ampere reflection}, \eqref{eq:pde for one source far field},
\eqref{eq:eq1} and \eqref{eq:pde for point source refraction} can be written as 
\begin{align}\label{eq:unique_eq}
\frac{d^{3/2}_{Q}(x,y)}{\sqrt{1-\phi_x^2(x,y)-\phi_y^2(x,y)}}|\det(D^2\phi(x,y))|=\frac{\tilde{f}\left(x,y\right)}{g\left(T(x,y)\right)},
\end{align}
where $\tilde{f}(x,y)=f(x,y)$ in the collimated case, and 
$\tilde{f}(x,y)=f\left(\frac{1}{\sqrt{x^2+y^2+1}}(x,y,1)\right)$ in the point source case.
\end{remark}
\section{Conclusion}
We have derived the equations that a phase discontinuity defined on a plane must satisfy in order that the resulting metasurface reflects or refracts light emanating with certain given variable intensity into a set of directions having also another variable intensity. 
The cases considered are when light emanates in a collimated beam and when light emanates from one point source. The resulting equations for the phase discontinuity are Monge-Amp\`ere partial differential equations which we show have solutions.

\end{document}